\documentclass[12pt]{article} 

\usepackage[utf8]{inputenc} 
\usepackage[margin=2.54cm]{geometry}
\usepackage[centertags]{amsmath}
\usepackage{amsfonts}
\usepackage{newlfont}
\usepackage{amsthm,amsmath,amssymb,amscd,verbatim,epsfig} 
\usepackage{makeidx}
\usepackage{layout}
\usepackage{color}

\theoremstyle{definition}
\newtheorem{thm}{Theorem}[section]
\newtheorem{cor}{Corollary}[section]
\newtheorem{lem}{Lemma}[section]

\newtheorem{ex}{Example}[section]

\theoremstyle{definition}
\newtheorem{defn}{Definition}[section]
\theoremstyle{remark}
\newtheorem{rem}{Remark}[section]

\numberwithin{equation}{section}

\begin{document}

\title{On the Hurwitz Zeta Function and Its Applications to Hyperbolic Probability Distributions}

\author{Tsung-Lin Cheng\footnote{Correspondent author: tlcheng@cc.ncue.edu.tw;
 National Changhua University of Education, Taiwan},\quad Chin-Yuan Hu\footnote{Retired Professor, National Changhua University of Education, Taiwan} }

\maketitle
\begin{abstract}
In this paper, we propose a new proof of the Jensen's formula (1895). We also derive some formulas similar to those in Pitman and Yor (2003). Besides, a new formula of the generalized Bernoulli function is also derived. At the end of the paper, the probability density functions of sinh- and tanh- are studied briefly for general cases.
\\

{\bf\emph{Keyword: Zeta function; Hyperbolic function; Ramanujan's Master Theorem\\
AMS Classification: 60E10}}
\end{abstract}
\newpage
\section{Introduction}
We can find a class of probability distributions generated by the Hurwitz Zeta function $\zeta(s,a)$ which is defined by the series
\begin{equation}
\zeta(s,a)=\sum\limits_{n=0}^{\infty}\frac1{(n+a)^s},
\end{equation}
where $s\in\mathbb{C}$, $Re(s)>1$, and $a>0$ is a real number (Apostol 1976). This function can be defined for all complex numbers except for the point $s=1$, at which it has a simple pole with residue 1, i.e.
\begin{equation}
  \lim\limits_{s\to 1}(s-1)\zeta(s,a)=1,\quad a>0.
 \end{equation}
  It also has the integral representation (Theorem 12.2 of Apostol 1976)
  \begin{equation}
    \Gamma(s)\zeta(s,a)=\int\limits_{0}^{\infty} \frac{x^{s-1}e^{-ax}}{1-e^{-x}}dx,\quad Re(s)>1,
  \end{equation}
where $\Gamma(\cdot)$ is the gamma function. Note that as $a=1$, $\zeta(s)\equiv\zeta(s,1)$ is the Riemann's Zeta function $\zeta(s)=\sum\limits_{n=1}^{\infty}\frac1{n^s}$, which has been investigated by Lin and Hu (2001).

Jensen (1895) obtained the following remarkable formula for the Riemann Zeta function
\begin{equation}
  (s-1)\zeta(s)=2\pi\int\limits_{-\infty}^{\infty} \frac{(\frac12+ix)^{1-s}}{(e^{\pi x}+e^{-\pi x})^2}dx,\quad s\in\mathbb{C}.
\end{equation}
He stated that this formula can easily be demonstrated with the help of the Cauchy's Theorem. Recently Johansson and Blagouchine (2019) proved this formula by using the contour integration method. Note that how Jensen actually computed $(s-1)\zeta(s)$ is not clear (see, Luschny, 2020).

Remember that the generalized Bernoulli function is defined by
\begin{equation}\label{Bern}
  B(s,a)=2\pi\int_{-\infty}^{\infty}\frac{(a-\frac12+ix)^s}{(e^{\pi x}+e^{-\pi x})^2}dx\equiv E((a-\frac12+iX)^s),\quad s\in\mathbb{C}, a>\frac12,
\end{equation}
where the random variable $X$ has the hyperbolic distribution function, that is, the probability function (p.d.f.) of $X$ is given by
\begin{equation}\label{hyper}
  f_X(x)=\frac{\pi}{2 cosh^2(\pi x)},\quad x\in \Re
\end{equation}
and the corresponding characteristic function is
\begin{equation}\label{hyperch}
  \phi_X(\theta)\equiv E(e^{i\theta X})=\frac{\theta/2}{sinh(\theta/2)},\quad \theta\in\Re,
\end{equation}
(see Pitman and Yor 2003). In particular, as $a=1$, we have
\[B(s)\equiv B(s,1)=2\pi \int\limits_{-\infty}^{\infty}\frac{(\frac12+ix)^s}{(e^{\pi x}+e^{-\pi x})^2}dx= E(\frac12+iX)^s,\quad s\in\mathbb{C}.\]
The Jensen's formula is the Bernoulli function  $B(1-s)=(s-1)\zeta(s)$, $s\in\mathbb{C}$. Note also that, the following relation holds:
$B(s,a)=-s\zeta(1-s,a),\quad s\in\mathbb{C},\quad a>1/2$ (see Theorem 12.13, Apostol 1976 ).

Pitman and Yor (2003) proposed three types of the hyperbolic random variables $\hat{S}_t$, $\hat{C}_t$, and $\hat{T}_t$. Their characteristic functions (c.f.) are given by
\begin{eqnarray*}
  \phi_{\hat{S}_t}(\theta)&\equiv & E(e^{i\theta \hat{S}_t})=\Big(\frac{\theta}{\sinh\theta}\Big)^t \\
   \phi_{\hat{C}_t}(\theta)&\equiv & E(e^{i\theta \hat{C}_t})=\Big(\frac{1}{\cosh\theta}\Big)^t \\
   \phi_{\hat{T}_t}(\theta)&\equiv & E(e^{i\theta \hat{T}_t})=\Big(\frac{\tanh\theta}{\theta}\Big)^t,\quad \theta\in\Re,\quad t>0.
\end{eqnarray*}
First, we note that $\hat{S}_t,\hat{C}_t$, and $\hat{T}_t$  are L{\'e}vy processes with finite moments of all orders. In particular, the distribution functions of
$\hat{S}_t,\hat{C}_t$, and $\hat{T}_t$ are of significance in analytic number theory for $t=1$ or $2$. For example, let $X$ be a random variable with the same distribution function as $\hat{S}_1/2$. Then the following identity holds,
\[B(s,a)=E(a-\frac12+iX)^s,\quad\hbox{ for } s\in\Re\hbox{ and }a>1/2.\]
It reveals that the generalized Bernoulli function is equal to a moment function of the hyperbolic r.v. $\hat{S}_1/2$. In particular, set $s=2n$ and $a=1$, we have
\[B_{2n}\equiv B(2n,1)=E(\frac12+iX)^{2n},\quad n=0,1,2,\cdots,\]
which implies the fact that the moment polynomial of $\hat{S}_1/2$ produces the Bernoulli number $B_{2n}$.
\begin{defn}
{\sl Let $X$ have the distribution function $F$ on $\Re$ and $a>0$ be a fixed real number. The Hurwitz Zeta function induced by the random variable $X$ is a complex valued function
defined by
\begin{equation}\label{HurZ}
  h_X(s,a)\equiv E(a+iX)^{-s}=\int\limits_{-\infty}^{\infty} \frac1{(a+ix)^s}dF(x),\quad s\ge 0.
\end{equation}
}
\end{defn}
The function $h_X$ is well-defined since the integral in the definition converges absolutely for $s\ge 0$, which we will prove in the following Lemma. Furthermore, if the
distribution function $F$ is symmetric, then the Hurwitz Zeta function can be written as
\[h_X(s,a)=2\int\limits_0^{\infty}\frac{\cos(s\tan^{-1}(x/a))}{(a^2+x^2)^{s/2}}dF(x),\quad s\ge 0, \quad a>0.\]
Note that if we impose the condition $E|X|^n<\infty$, for any $n\ge 1$, then the domain of $s$ can be extended to the complex plane $\mathbb{C}$.
Set $\theta=\tan^{-1}(x/a)$. We may write $a+ix=\sqrt{a^2+x^2}e^{i\theta}$.
\begin{eqnarray*}
  h_X(s,a)&=& \int\limits_{-\infty}^{\infty} (a+ix)^{-s}dF(x)=\int\limits_{-\infty}^{\infty} \Big[e^{i\theta}\sqrt{a^2+x^2}\Big]^{-s}dF(x) \\
  &=&\int\limits_{-\infty}^{\infty} \frac{[e^{-i\theta s}+e^{i\theta s}]/2}{(a^2+x^2)^{s/2}}dF(x)=\int\limits_{-\infty}^{\infty} \frac{\cos(s\theta)}{(a^2+x^2)^{s/2}}dF(x) \\
   &=& 2\int\limits_0^{\infty}\frac{\cos(s\tan^{-1}(x/a))}{(a^2+x^2)^{s/2}}dF(x).
\end{eqnarray*}
The representation is interesting and we can prove in the next section that, the Hurwitz Zeta function $h_X$ is an analytic function of $s\in\mathbb{C}$ for $a>0$.
This paper is organized as below: In Section 2, we will give the main results including a new proof of the Jensen's formula, the generalized Bernoulli function formula, and some new formulas related to the Bernoulli function. A class of new probability density functions on $(0,\infty)$ is also derived (see, Corollary 2.3), which is associated with the Hurwitz zeta function. In section 3, we will propose a new proof to modify the Master Theorem proposed by S. Ramanujan. Besides, we will give some interesting examples to illustrate our core ideas concerning the p.d.f.s of the hyperbolic cosh, sinh and tanh distribution mentioned in Pitman and Yor (2003). Our key approach is the Fourier inversion formula instead of the Cauchy residue theorem.

\section{Main Results}

\begin{lem}\label{ana}
{\sl Let $X$ have the distribution function $F$ on $\Re$ with all moments finite. Then the Hurwitz Zeta function $h_X$ is an analytic function of $s\in\mathbb{C}$ for $a>0$.}
\end{lem}
\begin{proof}
First, we verify the lemma for the case as $s\in\Re$, and then extend the result to complex $s$ by analytic continuation.
For $s\in\Re$ and $a>0$,
\[\Big|\int\limits_{-\infty}^{\infty}(a+ix)^{-s}dF(x)\Big|\le \int\limits_{-\infty}^{\infty}|a+ix|^{-s}dF(x)= \int\limits_{-\infty}^{\infty}(a^2+x^2)^{-\frac{s}2}dF(x). \]
If $s\ge 0$, then
\[\int\limits_{-\infty}^{\infty}(a^2+x^2)^{-s/2}dF(x)\le \int\limits_{-\infty}^{\infty}a^{-s}=a^{-s}<\infty. \]
On the other hand,if $s<0$, then $-s>0$, we may choose $n_0>-s>0$ and
\begin{eqnarray*}
  \int\limits_{-\infty}^{\infty}(a^2+x^2)^{-s/2}dF(x) &\le& \int\limits_{-\infty}^{\infty} [2\max\{a^2,x^2\}]^{-s/2}dF(x) \\
   &=& \int\limits_{a^2\le x^2} (2x^2)^{-s/2}dF(x)+\int\limits_{a^2> x^2} (2a^2)^{-s/2}dF(x)\\
   &\le& 2^{n_0/2}E(|X|^{-s})+2^{n_0/2}a^{-s}<\infty.
\end{eqnarray*}
Hence, the integral above converges absolutely for $s\in \Re$. Next, we show that the integral converges uniformly in every rectangle
\[\{s\in\mathbb{C}: s=\sigma+it,|\sigma|\le b, |t|\le c\},\quad b>0,c>0.\]
Observe that
\begin{eqnarray*}
&&\int\limits_{-\infty}^{\infty}|(a+ix)^{-s}|dF(x)=\int\limits_{-\infty}^{\infty}\Big|(\sqrt{a^2+x^2}e^{i\theta})^{-s}\Big|dF(x)\\
&=&\int\limits_{-\infty}^{\infty}\Big|(\sqrt{a^2+x^2})^{-\sigma-it}e^{i\theta (-\sigma-it)}\Big|dF(x)\\
&=&\int\limits_{-\infty}^{\infty}(a^2+x^2)^{-\sigma/2}e^{t\theta}dF(x)
\le e^{\frac{\pi}{2}c}\int\limits_{-\infty}^{\infty}(a^2+x^2)^{-\sigma/2}dF(x).
\end{eqnarray*}
If $\sigma>0$, $\sigma\le b$, then
\begin{eqnarray*}
 && \int\limits_{-\infty}^{\infty}|(a+ix)^{-s}|dF(x) \le e^{c\pi/2}\int\limits_{-\infty}^{\infty}(a^2+x^2)^{-\sigma/2}dF(x) \\
   &\le&   e^{c\pi/2}\int\limits_{-\infty}^{\infty}a^{-\sigma}dF(x)\le e^{c\pi/2} \Big(max\{1,1/a\}\Big)^b.
\end{eqnarray*}
If $\sigma<0$, $|\sigma|\le b$, then we may choose a positive integer $n_0>b$ and
\begin{eqnarray*}
   && \int\limits_{-\infty}^{\infty}|(a+ix)^{-s}|dF(x)\le e^{c\pi/2}\int\limits_{-\infty}^{\infty}(a^2+x^2)^{-\sigma/2}dF(x) \\
   &\le& e^{c\pi/2}\int\limits_{-\infty}^{\infty} [2\max\{a^2,x^2\}]^{-\sigma/2}dF(x)\\
   &=& e^{c\pi/2}\Big[\int\limits_{\{a^2\le x^2\}}[2\max\{a^2,x^2\}]^{-\sigma/2}dF(x)+\int\limits_{\{a^2> x^2\}}[2\max\{a^2,x^2\}]^{-\sigma/2}dF(x)\Big] \\
   &\le& e^{c\pi/2}\Big[\int\limits_{-\infty}^{\infty}2^{-\sigma/2}(x^2)^{-\sigma/2}dF(x)+ \int\limits_{-\infty}^{\infty}(2a^2)^{-\sigma/2}dF(x)\Big]\\
   &\le& e^{c\pi/2} [2^{n_0/2}E(|X|^{-\sigma})+2^{n_0/2}a^{-\sigma}]\\
   &\le& e^{c\pi/2}2^{n_0/2}\Big\{\max\{1,E(|X|^{n_0})\}+(\max\{1,a\})^b\Big\}.
\end{eqnarray*}
The above argument shows that the integral converges uniformly in every rectangle $\{s\in\mathbb{C}:s=\sigma+it, |\sigma|\le b, |t|\le c\}$, for $b>0$ and $c>0$.
By Analytic Continuation Theorem, we can prove that the integral is analytic for all $s\in\mathbb{C}$.
\end{proof}
Next, we will introduce three type of hyperbolic moment functions by six series as below.
For $a>0$, we define
\begin{itemize}
  \item[(1)]
  \begin{description}
    \item[] $S_1(s,a)=\sum\limits_{n=0}^{\infty} \frac1{(n+a)^s}$, $s>1$.
    \item[] $S_2(s,a)=\sum\limits_{n=0}^{\infty} \frac{n}{(n+a)^s}$, $s>2$.
  \end{description}
 \item[(2)]
  \begin{description}
  \item[] $C_1(s,a)=\sum\limits_{n=0}^{\infty} \frac{(-1)^n}{(n+a)^s}$, $s>1$.
 \item[] $C_2(s,a)=\sum\limits_{n=0}^{\infty} \frac{(-1)^n(n+1)}{(n+a)^s}$, $s>2$.
\end{description}
 \item[(3)]
  \begin{description}
    \item[] $T_1(s,a)=\sum\limits_{n=0}^{\infty} 2(-1)^n E(a+U+n)^{-s}$, $s>1$, and
   \item[] $T_2(s,a)=\sum\limits_{n=1}^{\infty} 4(-1)^{n-1}n E(a-1+V+n)^{-s}$, $s>2$,

   where $U$ has the uniform distribution over $(0,1)$ and $V$ has the triangular distribution over $(0,2)$, namely, the pdf of $V$ is
    \begin{equation*}
    f_V(u)=\left\{\begin{array}{ll}
                 u,\quad & \mbox{if}\quad 0\le u\le 1  \\  
                 2-u,\quad &\mbox{if}\quad 1<u\le 2\\
                 0,\quad & \mbox{ otherwise}.   
                               \end{array} \right.
    \end{equation*}

  \end{description}
\end{itemize}
The following theorem is a restatement of Jensen (1895) and Johansson and Blagouchine (2019), for which we give a new derivation by using both Lemma 2.1 and Fourier transform.
\begin{thm}
{\sl The Hurwitz Zeta function $\zeta(s,a)$ and the hyperbolic sinh-moment function $S_1(s,a)$ can be related by the following equation
\begin{equation}\label{JJB1}
  (s-1)\zeta(s,a) = (s-1)S_1(s,a)=E(a-\frac12+iX)^{1-s},\quad\hbox{for }a>1/2, s\in\mathbb{C}.
 \end{equation}
In particular, when $a=1$, we have
\begin{equation}\label{JJB2}
  (s-1)\zeta(s)=(s-1)S_1(s,1)=E(\frac12+iX)^{1-s},\quad s\in\mathbb{C},
\end{equation}
where the random variable $X$ has the hyperbolic sinh-distribution as $\hat{S}_1/2$, i.e., the pdf
of $X$ is $f_X(x)=\frac{\pi}{2\cosh^2(\pi x)}$, $x\in\Re$.
}
\end{thm}
\begin{proof}
Observe that for $s>1$, $a>1/2$, the Hurwitz Zeta function has the integral representation

\begin{eqnarray*}
  &&\Gamma(s)\zeta(s,a) = \int\limits_0^{\infty} x^{s-1}\frac{e^{-ax}}{1-e^{-x}}dx \\
  &=& \int\limits_0^{\infty} x^{s-2}e^{-(a-\frac12)x}\frac{xe^{-x/2}}{1-e^{-x}}dx\\
  &=& \int\limits_0^{\infty} x^{s-2}e^{-(a-\frac12)x} \frac{x/2}{\sinh{\frac{x}2}}dx\\
  \end{eqnarray*}
\begin{eqnarray*}
 &=& \int\limits_0^{\infty} x^{s-2}e^{-(a-\frac12)x}\int\limits_{-\infty}^{\infty} e^{itx}dF_X(t)dx,\quad\hbox{( by \ref{hyperch})}\\
    &=&\int\limits_{-\infty}^{\infty}\int\limits_0^{\infty}x^{s-2}e^{-(a-\frac12)x}e^{itx}dx\,dF_X(t) \\
  &=& \int\limits_{-\infty}^{\infty} \frac{\Gamma(s-1)}{(a-\frac12-it)^{s-1}}dF_X(t).
\end{eqnarray*}
Substituting $\Gamma(s)=(s-1)\Gamma(s-1)$ into the left-hand-side of the above equation and dividing both side by $\Gamma(s-1)$, we have
\begin{equation}
(s-1)\zeta(s,a)=\int\limits_{-\infty}^{\infty}(a-\frac12-it)^{1-s}dF_X(t)=E(a-\frac12-iX)^{1-s}=E(a-\frac12+iX)^{1-s}.
\end{equation}

Note that the random variable $X{\buildrel d \over{=}} \hat{S}_1/2$ has the moments of all orders, we have
\[(s-1)\zeta(s,a)=E(a-\frac12+iX)^{1-s}=h_X(s-1,a-\frac12),\quad s\in \mathbb{C},\]
where the function $f_X$ is the Hurwitz Zeta function induced by the random variable $X$. By Lemma \ref{ana}, the function represents an analytic
function for $s\in \mathbb{C}$.
\end{proof}
Now we immediately have the following corollary just by interchanging $s$ with $1-s$ in the equation (\ref{JJB1}).
\begin{cor}
{\sl Let $X{\buildrel d \over{=}} \hat{S}_1/2$ be a hyperbolic $\sinh$-distributed random variable with pdf
(\ref{hyper}). Then the generalized Bernoulli function is
\begin{equation}\label{gBern}
B(s,a)=E(a-\frac12+iX)^s,\quad s\in \mathbb{C}, a>\frac12.
\end{equation}
In particular, when $a=1$, we have the Bernoulli function $B(s)=E(\frac12+iX)^s, \quad s\in \mathbb{C}$.
 }
\end{cor}
It is interesting that $2\pi X {\buildrel d \over{=}}\pi\hat{S}_1$ has the logistic distribution
\[P(2\pi X>x)=\frac1{1+e^{-x}},\quad x\in\Re,\]
which has found numerous applications. The following theorem is a new integral representation of the Hurwitz Zeta function.
\begin{thm}\label{HurZ2}
{\sl Let $U_1$ be a random variable with the pdf $f_1(u)=1$, $0<u<1$, and zero otherwise, and $U_2$ be another random variable with pdf
$f_2(u)=2(1-u)$, $0<u<1$, and zero otherwise. Suppose that $X{\buildrel d \over{=}} \hat{S}_2/2 $ has the hyperbolic sinh-distribution with the pdf
\[f_X(x)=\frac{\pi [\pi x \coth(\pi x)-1]}{\sinh^2(\pi x)},\quad x\in\Re.\]
Then the Hurwitz Zeta function $\zeta(s,a)$ has the integral representation
\begin{equation}\label{Hurz3}
(s-1)\zeta(s,a)=\frac12 E(a-1+U_2+iX)^{-s}+(a-1)E(U_1+iX)^{-s},\quad\hbox{ for }a\ge 1, s\in\mathbb{C}.
\end{equation}
In particular, when $a=1$, we have
\begin{equation}\label{HurZ4}
(s-1)\zeta(s)=\frac12 E(U_2+iX)^{-s},\quad s\in\mathbb{C}.
\end{equation}
}
\end{thm}
\begin{proof}
The Hurwitz Zeta function has the integral representation
\[\Gamma(s)\zeta(s,a) = \int\limits_0^{\infty} x^{s-1}\frac{e^{-ax}}{1-e^{-x}}dx,\quad s>1,a>0.\]
By integration by parts, we have
\[\Gamma(s)\zeta(s,a)=\frac1{s} \int\limits_0^{\infty} x^{s}\frac{e^{-ax}[a-(a-1)e^{-x}]}{(1-e^{-x})^2}dx, \]
 which is
 \[s\Gamma(s)\zeta(s,a)=\int\limits_0^{\infty} x^{s}\frac{e^{-ax}[a-(a-1)e^{-x}]}{(1-e^{-x})^2}dx. \]
 Finally, by subtraction we have
 \begin{eqnarray*}
    (s-1)\Gamma(s)\zeta(s,a)&=&\int\limits_0^{\infty}x^{s-1}e^{-ax}\frac{[ax-1+e^{-x}(1-(a-1)x)]}{(1-e^{-x})^2}dx\\
    &=& \int\limits_0^{\infty}x^{s-1}e^{-(a-1)x}\Big[\frac12\frac{e^{-x}-1+x}{x^2/2}+(a-1)\frac{1-e^{-x}}{x}\Big]\Big(\frac{x/2}{\sinh{\frac{x}2}}\Big)^2dx,
 \end{eqnarray*}
 which has finite values as $s>1$ and $a\ge 1$.
 Write $(s-1)\Gamma(s)\zeta(s,a)=I_1+I_2$, where
 \begin{eqnarray*}
 I_1&=&\int\limits_0^{\infty}x^{s-1}e^{-(a-1)x}\Big[\frac12\frac{e^{-x}-1+x}{x^2/2}\Big]\Big(\frac{x/2}{\sinh{\frac{x}2}}\Big)^2dx,\quad\hbox{ and}\\
 I_2&=&\int\limits_0^{\infty}x^{s-1}e^{-(a-1)x}\Big[(a-1)\frac{1-e^{-x}}{x}\Big]\Big(\frac{x/2}{\sinh{\frac{x}2}}\Big)^2dx.
 \end{eqnarray*}
 Now the moment generation functions of $U_1$ and $U_2$, respectively, and the characteristic function of $X$ are given by
 \begin{eqnarray*}
   E(e^{-\lambda U_1}) &\equiv& \int_{\Re} e^{-\lambda u}dF_{U_1}(u)=\frac{1-e^{-\lambda}}{\lambda},\quad\lambda>0,  \\
     E(e^{-\lambda U_2}) &\equiv& \int_{\Re} e^{-\lambda u}dF_{U_2}(u)=\frac{e^{-\lambda}-1+\lambda}{\lambda^2/2},\quad\lambda>0,\\
     E(e^{i\theta X})&\equiv& \int_{\Re} e^{i\theta x}dF_X(x)=\Big(\frac{\theta/2}{\sinh{\frac{\theta}2}}\Big)^2,\quad\theta\in\Re .
 \end{eqnarray*}
 Note that
 \begin{eqnarray*}
 I_1&=&\int\limits_0^{\infty}x^{s-1}e^{-(a-1)x}\frac12 \int_{\Re} e^{-\lambda u}dF_{U_2}(u)\int_{\Re} e^{it x}dF_X(t)\,dx\\
 &=& \frac12\int_{\Re}\int_{\Re}\int\limits_0^{\infty}x^{s-1}e^{-[(a-1)+u-it]x}dx\,dF_{U_2}(u)\,dF_X(t)\\
 &=& \frac12\int_{\Re}\int_{\Re}\frac{\Gamma(s)}{[(a-1)+u-it]^s}dF_{U_2}(u)\,dF_X(t)\\
 &=&  \frac12\Gamma(s)E(a-1+U_2-iX)^{-s}\\
 &=& \frac12\Gamma(s)E(a-1+U_2+iX)^{-s},\quad \hbox{ by symmetry of the distribution of }X.
 \end{eqnarray*}
 In a similar fashion, we may obtain
 \[I_2=\Gamma(s)(a-1)E(U_1+iX)^{-s}.\]
 Therefore,
 \begin{eqnarray*}
   (s-1)\zeta(s,a) &=& \frac12 E(a-1+U_2+iX)^{-s}+(a-1)E(U_1+iX)^{-s}.
 \end{eqnarray*}
 By analytic continuation, this equation holds for all $s\in\mathbb{C}$ and $a\ge 1$.
\end{proof}
Observe that,
\begin{eqnarray*}
  B(s,a) &=& -s\zeta(1-s,a)\hbox{ and }
  B(s) \equiv B(s,1)=-s\zeta(1-s),\quad\hbox{ or }  \\
   B(1-s,a)&=& (s-1)\zeta(s,a).
\end{eqnarray*}
We immediately have the following corollary.
\begin{cor}
{\sl Let $U_1,U_2$ and $X$ be the random variables as given in Theorem \ref{HurZ2}. The the generalized Bernoulli function $B(s,a)$
has the integral representation, for $a\ge 1$, $s\in\mathbb{C}$,
\begin{equation}
B(s,a)=\frac12 E(a-1+U_2+iX)^{s-1}+(a-1)E(U_1+iX)^{s-1}.
\end{equation}
In particular, when $a=1$, we have the Bernoulli function
\begin{equation}\label{Bf}
  B(s)=\frac12 E(U_2+iX)^{s-1}, \quad s\in\mathbb{C}.
\end{equation}
}
\end{cor}
From the proof of Theorem \ref{HurZ2}, for $a>0$, we have $\lim\limits_{s\to 1}(s-1)\zeta(s,a)=1$, which implies that
$h_a(x)$ is a p.d.f. and we have the following corollary.
\begin{cor}
{\sl The following identity holds for $a\ge 1$ and $Re(s)>0$
\begin{equation}\label{idhurz}
  \int\limits_0^{\infty} x^{s-1}h_a(x)dx=(s-1)\Gamma(s)\zeta(s,a),
\end{equation}
where $h_a(x)$ is a p.d.f. on $(0,\infty)$ defined by
\[h_a(x)=\frac{e^{-ax}\Big[ax-1+e^{-x}\Big(1-(a-1)x\Big)\Big]}{(1-e^{-x})^2},\quad x>0.\]
}
\end{cor}
\begin{rem}
  The p.d.f. $h_a(x)$ is characterized by the above Mellin transform where the right side is interpreted by continuity at $s=1$. Aldous (2001) obtained  results for the
  the special case as $a=1$, $s=1$ and $s=2$. Pitman and Yor (2003) studied the special case at $a=1$. In Pitman and Yor (2003), the following identity holds for $x>0$,
  \[h_1(x)=P(\pi \hat{S}_2 >x)=\frac{e^{-x}(e^{-x}-1+x)}{(1-e^{-x})^2}.\]
\end{rem}
\begin{thm}\label{HurZ5}
  {\sl Let $X{\buildrel d\over{=}}\hat{S}_2/2$ be the hyperbolic sinh-distributed random variable with the p.d.f.
  \begin{equation}\label{hs}
    f_X(x)=\frac{\pi [\pi x \coth(\pi x)-1]}{\sinh^2(\pi x)},\quad x\in\Re.
  \end{equation}
  Then the hyperbolic sinh-moment function $S_2(s,a)$ has the integral representation
  \begin{equation}\label{hsm}
    S_2(s,a)=\frac1{(s-1)(s-2)}E(a+iX)^{2-s},\quad Re(s)\ne 1 \hbox{ or } 2,\quad a>0.
  \end{equation}
   }
\end{thm}
\begin{proof}
The hyperbolic sinh-moment function $S_2(s,a)$ is defined by
\[S_2(s,a)=\sum\limits_{n=0}^{\infty} \frac{n}{(n+a)^s},\quad s>2,\quad a>0.\]
Note that since this series is absolutely convergent for $s$ with $Re(s)>2$, we can interchange the sum and integral in the following
calculation
\begin{eqnarray*}
   && \Gamma(s)S_2(s,a)=\sum\limits_{n=0}^{\infty} \frac{n\Gamma(s)}{(n+a)^s} \\
   &=&\sum\limits_{n=1}^{\infty} n\int\limits_0^{\infty} x^{s-1}e^{-(n+a)x}dx
   = \int\limits_0^{\infty} x^{s-1}\sum\limits_{n=1}^{\infty}n e^{-nx}\cdot e^{-ax}dx \\
   &=& \int\limits_0^{\infty} x^{s-1}e^{-ax}\frac{e^{-x}}{(1-e^{-x})^2}dx
   =  \int\limits_0^{\infty} x^{s-3}e^{-ax}\Big(\frac{x/2}{\sinh{x/2}}\Big)^2dx
   \end{eqnarray*}
   \begin{eqnarray*}
   &=&  \int\limits_0^{\infty} x^{s-3}e^{-ax}\int\limits_{\Re}e^{ixt}dF_X(t)dx=\int\limits_{\Re}\int\limits_0^{\infty} x^{s-3}e^{-(a-it)x}dx\,dF_X(t)\\
   &=& \int\limits_{\Re} \frac{\Gamma(s-2)}{(a-it)^{s-2}}dF_X(t)=\Gamma(s-2)E(a-iX)^{2-s}\\
   &=&\Gamma(s-2)E(a+iX)^{2-s},\quad\hbox{by symmetry of the distribution of $X$}.
\end{eqnarray*}
Furthermore, we may obtain the equality $(s-1)(s-2)S_2(s,a)=E(a+iX)^{2-s}$. By analytic continuation, we can extend the above equation to $s\in\mathbb{C}$, $a>0$.
\end{proof}
The next theorem is concerning the integral representation of the hyperbolic cosh-moment function.
\begin{thm}\label{hycosh}
  {\sl Let $X{\buildrel d \over{=}}\hat{C}_1/2$ be the hyperbolic cosh-distributed random variable with the p.d.f.
  \[f_X(x)=\frac1{\cosh(\pi x)},\quad x\in \Re.\]
  Then the hyperbolic cosh-moment function $C_1(s,a)$ has the integral representation
  \begin{equation}\label{hycoshm}
    C_1(s,a)=\frac12 E(a-\frac12+iX)^{-s},\quad\hbox{for }a>\frac12\hbox{ and }s\in\mathbb{C}.
  \end{equation}
  In particular, as $a=1$, we have
  \[(1-2^{1-s})\zeta(s)=C_1(s,1)=\frac12 E(\frac12+iX)^{-s},\quad s\in\mathbb{C}.\]
    }
\end{thm}
\begin{proof}
Note that $C_1(s,a)=\sum\limits_{n=0}^{\infty} \frac{(-1)^n}{(n+a)^s}$, $s>1$, $a>\frac12$, and we have
\begin{eqnarray*}
   && \Gamma(s)C_1(s,a)=\sum\limits_{n=0}^{\infty}\frac{(-1)^n \Gamma(s)}{(n+a)^s}=\sum\limits_{n=0}^{\infty}(-1)^n\int\limits_{0}^{\infty} x^{s-1}e^{-(n+a)x}dx \\
    &=&\sum\limits_{n=0}^{\infty} (-1)^n\int\limits_{0}^{\infty} x^{s-1}e^{-(n+a)x}dx
   = \int\limits_{0}^{\infty} x^{s-1}e^{-ax}\sum\limits_{n=0}^{\infty}(-e^{-x})^ndx\\
   &=&  \int\limits_{0}^{\infty} x^{s-1}e^{-ax}\frac1{1+e^{-x}}dx= \int\limits_{0}^{\infty} x^{s-1}e^{-(a-\frac12)x} \frac1{2\cosh(\frac{x}2)}dx\\
   &=& \frac12\int\limits_{0}^{\infty} x^{s-1}e^{-(a-\frac12)x}\int\limits_{\Re}e^{itx}dF_X(t)dx,\quad\hbox{ by the definition of }X,\\
   &=&\frac12\int\limits_{\Re}\int\limits_{0}^{\infty}x^{s-1}e^{-(a-\frac12-it)x}dx\,dF_X(t)
   =\frac12\int\limits_{\Re}\frac{\Gamma(s)}{(a-\frac12-it)^s}dF_X(t)\\
   &=&\frac{\Gamma(s)}2 E(a-\frac12+iX)^{-s},\quad\hbox{by symmetry of the distribution of }X.
   \end{eqnarray*}
Dividing both sides by $\Gamma(s)$ and by analytic continuation, we have
\[C_1(s,a)=\frac12E(a-\frac12+iX)^{-s},\quad\hbox{for }a>1/2, s\in\mathbb{C}.\]
\end{proof}
The following corollary is an interesting result of the above theorem. By the defintion of the series $C_1(s,a)$, we have, as $a=1$,
\begin{equation*}
C_1(s,1)=\sum\limits_{n=1}^{\infty} \frac{(-1)^{n-1}}{n^s}=(1-2^{1-s})\zeta(s),\quad Re(s)>0
\end{equation*}
(see, Apostol 1976, p292), while as $a=\frac12$,
\begin{equation*}
C_1(s,\frac12)=2^s\sum\limits_{n=0}^{\infty}\frac{(-1)^n}{(2n+1)^s}\equiv 2^sL_{\chi_4}(s),\quad Re(s)>0,
\end{equation*}
where $L_{\chi_4}(s)$ is the Dirichlet series associated with the quadratic character modulo 4 (see, Biane, Pitman and Yor 2001).
\begin{cor}
  {\sl The following identity holds for $a>1/2$ and $Re(s)>-1$,
  \[2\Gamma(s+1)C_1(s,a)=\int\limits_{\Re}|x|^sg_a(x)dx,\]
  where $g_a(x)$ is a p.d.f. $g_a(x)=\frac{e^{-a|x|}(a+(a-1)e^{-|x|})}{(1+e^{-|x|})^2}$, $x\in \Re$.

  In particular, when $a=1$, we have the absolute moment formula of the logistic distribution, namely,
  \begin{equation}\label{logist}
    E(|X|^s)=2\Gamma(s+1)(1-2^{1-s})\zeta(s),\quad \Re(s)>-1,
  \end{equation}
  where $X$ has the logistic distribution $F_X(x)=\frac1{1+e^{-x}},\quad x\in\Re$.
  }
\end{cor}
\begin{proof}
By Theorem \ref{hycosh} and integration by parts, we have
\begin{eqnarray*}
  \Gamma(s)C_1(s,a) &=& \frac1{s}\int_0^{\infty}\frac{e^{-ax}}{1+e^{-x}}dx^s \\
   &=& \frac1{s} \int_0^{\infty} x^s\frac{e^{-ax}[a+(a-1)e^{-x}]}{(1+e^{-x})^2}dx,\quad\hbox{ which is equivalent to}\\
   \Gamma(s+1)C_1(s,a)&=&\int_0^{\infty}x^s\frac{e^{-ax}[a+(a-1)e^{-x}]}{(1+e^{-x})^2}dx.
\end{eqnarray*}
Since $\lim\limits_{s\to 0}C_1(s,a)=1/2$, we have
\[2\Gamma(s+1)C_1(s,a)=\int\limits_{\Re}|x|^s g_a(x)dx,\quad Re(s)>-1.\]
Furthermore, as $g_a(x)\ge 0$, for $x\in\Re$, and $a>1/2$, we can see that $g_a(x)$ is a p.d.f. The second assertion of the corollary can be
treated as a special case of the above theorem as $a=1$.
\end{proof}
\begin{thm}\label{HurZ6}
{\sl Let $X\buildrel d\over{=}\hat{C}_2/2$ be a random variable with the hyperbolic cosh-distribution, i.e. its p.d.f. can be written as
\[f_X(x)=\frac{2x}{\sinh(\pi x)},\quad x\in \Re .\] Then, the hyperbolic cosh-moment function $C_2(s,a)$ has the integral representation
\begin{equation}\label{hycoshm}
  C_2(s,a)=\frac14E(a-1+iX)^{-s},\quad s\in\mathbb{C},\quad a>1.
\end{equation}
}
\end{thm}
\begin{proof}
  By the definition of the hyperbolic cosh-moment function $C_2(s,a)=\sum\limits_{n=0}^{\infty}\frac{(-1)^n(n+1)}{(n+a)^s}$, $s>2$, we have
  \begin{eqnarray*}
   && \Gamma(s)C_2(s,a) = \sum\limits_{n=0}^{\infty} \frac{(-1)^n(n+1)\Gamma(s)}{(n+a)^s} \\
     &=& \sum\limits_{n=0}^{\infty}\binom{-2}{n}\int\limits_{0}^{\infty} x^{s-1}e^{-(n+a)x}dx=\int\limits_{0}^{\infty} x^{s-1}e^{-ax} \sum\limits_{n=0}^{\infty}\binom{-2}{n}
     e^{-nx}dx\\
  &=&\int\limits_{0}^{\infty} x^{s-1}e^{-ax}\Big(\frac1{1+e^{-x}}\Big)^2dx,\\
  &&\quad\hbox{ which is due to the equality }
    (1+x)^{-2}=\sum\limits_{n=0}^{\infty}\binom{-2}{n}x^n.
 \end{eqnarray*}
Besides, by the definition of $X$, we have
\[E(e^{i\theta X})=\int\limits_{\Re}e^{i\theta x}dF_X(x)=\Big(\frac1{\cosh{\frac{\theta}2}}\Big)^2,\quad\theta\in\Re,\]
which implies
\begin{eqnarray*}
\Gamma(s)C_2(s,a)&=&\int\limits_{0}^{\infty} x^{s-1}e^{-ax}\Big(\frac1{1+e^{-x}}\Big)^2dx\\
&=& \int\limits_{0}^{\infty} x^{s-1}e^{-(a-1)x}\frac14\Big(\frac1{\cosh\frac{x}2}\Big)^2dx\\
&=&\frac14 \int\limits_{0}^{\infty} x^{s-1}e^{-(a-1)x}\int_{\Re}e^{itx}dF_X(t)dx\\
&=&\frac14\int\limits_{\Re} \int\limits_{0}^{\infty}x^{s-1}e^{-(a-1-it)x}dxdF_X(t)\\
&=&\frac14\int\limits_{\Re} \frac{\Gamma(s)}{(a-1-it)^s}dF_X(t)=\frac14\Gamma(s)E(a-1-iX)^{-s}\\
&=&\frac14\Gamma(s)E(a-1+iX)^{-s}.
\end{eqnarray*}
Dividing both sides by $\Gamma(s)$ and using analytic continuation, we complete the proof.
\end{proof}
We will prove the following theorem concerning hyperbolic tanh-distribution function.
\begin{thm}\label{hurw7}
{\sl Let $X\buildrel d\over{=}\hat{T}_1/2$ have the hyperbolic tanh-distribution with the p.d.f.
\[f_X(x)=\frac2{\pi} \log\coth(\frac{\pi}2 |x|),\quad x\in\Re.\]
Then, the hyperbolic tanh-moment function $T_1(s,a)$ has the integral representation
\begin{equation}\label{tanhm}
T_1(s,a)=E(a+iX)^{-s},\quad\hbox{for }a>0,\quad s\in\mathbb{C}.
\end{equation}
}
\end{thm}
\begin{proof}
By the definition of the hyperbolic tanh-moment function $T_1(s,a)$, multiplying it by $\Gamma(s)$, we have
\begin{eqnarray*}
  &&\Gamma(s)T_1(s,a)=\Gamma(s)\sum\limits_{n=0}^{\infty} 2(-1)^nE(a+U+n)^{-s}\\
  &=&\sum\limits_{n=0}^{\infty} 2(-1)^n \int\limits_{\Re}\frac{\Gamma(s)}{(a+u+n)^s}dF_U(u),\quad\hbox{ in which } U\buildrel d\over{=}Unif(0,1)\\
  &=&\sum\limits_{n=0}^{\infty} 2(-1)^n \int\limits_{\Re}\int\limits_0^{\infty} x^{s-1}e^{-(a+u+n)x}dxdF_U(u)\\
  &=&\int\limits_0^{\infty}x^{s-1}e^{-ax}2\int\limits_{\Re}e^{-ux}dF_U(u)\sum\limits_{n=0}^{\infty} (-1)^ne^{-nx}dx\\
  &=&\int\limits_0^{\infty}x^{s-1}e^{-ax}2\frac{1-e^{-x}}{x}\cdot \frac1{1+e^{-x}}dx\\
  &=&\int\limits_0^{\infty}x^{s-1}e^{-ax}\frac{\tanh(\frac{x}2)}{x/2}dx.
\end{eqnarray*}
Note that $X\buildrel d\over{=}\hat{T}_1/2$. Therefore,
\begin{eqnarray*}
   && \int\limits_0^{\infty}x^{s-1}e^{-ax}\frac{\tanh(\frac{x}2)}{x/2}dx \\
   &=& \int\limits_0^{\infty}x^{s-1}e^{-ax}\int\limits_{\Re}e^{itx}dF_X(t)dx\\
   &=& \int\limits_{\Re}\int\limits_0^{\infty}x^{s-1}e^{-(a-it)x}dx\,dF_X(t)\\
   &=& \int\limits_{\Re} \frac{\Gamma(s)}{(a-it)^s}dF_X(t)=\Gamma(s)E(a+iX)^{-s}.
\end{eqnarray*}
Dividing both sides by $\Gamma(s)$ and by analytic continuation, we have
\[T_1(s,a)=E(a+iX)^{-s},\quad\hbox{for }a>0,\quad s\in\mathbb{C}.\]
\end{proof}
\begin{thm}
Let $X\buildrel d\over{=} \hat{T}_2/2$ have the hyperbolic tanh-distribution with the p.d.f.
\begin{equation}\label{hytanh}
  f_X(x)=\int\limits_{2|x|}^{\infty}\frac{y(y-2|x|)}{\sinh(\frac{\pi}2 y)}dy,\quad x\in\Re.
\end{equation}
Then the hyperbolic tanh-moment function $T_2(s,a)$ has the integral representation
\begin{equation}\label{hytanhm}
  T_2(s,a)=E(a+iX)^{-s},\quad a>0,s\in\mathbb{C}.
\end{equation}
\end{thm}
\begin{proof}
  The hyperbolic tanh-moment function $T_2(s,a)$ is defined by
  \[T_2(s,a)=\sum\limits_{n=1}^{\infty}4(-1)^{n-1}nE(a-1+V+n)^{-s},\quad s>2,\quad a>0,\]
  where $V$ has the triangle p.d.f., namely, the pdf of $V$ is
    \begin{equation*}
    f_V(u)=\left\{\begin{array}{ll}
                 u,\quad & \mbox{if}\quad 0\le u\le 1  \\  
                 2-u,\quad &\mbox{if}\quad 1<u\le 2\\
                 0,\quad & \mbox{ otherwise}.   
                               \end{array} \right.
    \end{equation*}
Observe that
\begin{eqnarray*}
   && \Gamma(s)T_2(s,a)=\sum\limits_{n=1}^{\infty}4(-1)^{n-1}n\Gamma(s)E(a-1+V+n)^{-s} \\
   &=& \int\limits_{\Re}\sum\limits_{n=1}^{\infty}4(-1)^{n-1}n\frac{\Gamma(s)}{[(a-1)+v+n]^s}dF_V(v)  \\
   &=& \int\limits_{\Re}\int\limits_0^{\infty} \sum\limits_{n=1}^{\infty}4(-1)^{n-1}nx^{s-1}e^{-(a-1+v+n)x}dxdF_V(v) \\
   &=& \int\limits_0^{\infty} x^{s-1}e^{-(a-1)x}\Big[\int\limits_{\Re}e^{-vx}dF_V(v)\Big] \sum\limits_{n=1}^{\infty}4(-1)^{n-1}ne^{-nx}dx\\
   &=& \int\limits_0^{\infty} x^{s-1}e^{-(a-1)x}\Big(\frac{1-e^{-x}}{x}\Big)^2 \frac{4e^{-x}}{(1+e^{-x})^2}dx\\
   &=& \int\limits_0^{\infty}x^{s-1}e^{-ax}\Big(\frac{\tanh{\frac{x}2}}{x/2}\Big)^2dx.
\end{eqnarray*}
By the definition of $X$, its characteristic function is given by
\[E(e^{i\theta X})=\int\limits_{\Re}e^{i\theta x}dF_X(x)=\Big(\frac{\tanh{\frac{\theta}2}}{\theta/2}\Big)^2,\quad\theta\in\Re,\]
which implies
\begin{eqnarray*}
   && \Gamma(s)T_2(s,a)=\int\limits_0^{\infty} x^{s-1}e^{-ax}\Big(\frac{\tanh{\frac{x}2}}{x/2}\Big)^2dx\\
   &=&=\int\limits_0^{\infty} x^{s-1}e^{-ax}\int\limits_{\Re}e^{it x}dF_X(t)\,dx \\
   &=& \int\limits_{\Re}\int\limits_0^{\infty} x^{s-1}e^{-(a-it)x}dx\,dF_X(t)\\
   &=& \int\limits_{\Re}\frac{\Gamma(s)}{(a-it)^s}dF_X(t)=\Gamma(s)E(a-iX)^{-s}=\Gamma(s))E(a+iX)^{-s},
\end{eqnarray*}
and $T_2(s,a)=E(a+iX)^{-s}$. Again by analytic continuation, we prove the theorem.
\end{proof}
\section{Some More General Results Concerning Ramanujan's Master Theorem and Pitman-Yor Hyperbolic Distributions}
S. Ramanujan introduced the so called Ramanujan's Master Theorem that provides an explicit expression for the
Mellin transform of a function. It was widely used as a tool in computing definite integrals and infinite series (see Berndt 1985). As stated by Ramanujan in his quarterly reports, his
Master Theorem is given by the following interesting relation
\begin{equation}\label{RamaMaster}
  \int\limits_0^{\infty} x^{s-1}\sum\limits_{j=0}^{\infty}\frac{(-x)^j\lambda(j)}{j!}dx=\Gamma(s)\lambda(-s).
\end{equation}
Another interesting identity by Ramanujan alongside (\ref{RamaMaster}) is given by
\begin{equation}\label{RamaMaster2}
  \int\limits_0^{\infty} x^{s-1}\sum\limits_{j=0}^{\infty}(-x)^j\phi(j)dx=\frac{\pi}{\sin(\pi s)}\phi(-s),
\end{equation}
where $\phi(n)=\lambda(n)/\Gamma(n+1)$.

Amdeberhan et. al. (2012) mentioned that, ``the proof of Ramanujan's Master Theorem provided by Hardy (1978)
employs Cauchy's residue theorem as well as the well-known Mellin inversion formula..." As Hardy observed, the Fourier inversion lay outside the range of Ramanujan's ideas (see Edwards 1974, p.223-224). Therefore it is necessary to modify the Original Ramanujan's Master Theorem. In the following, we propose a new proof to modify the Master Theorem by using the formula
of the characteristic function together with the Laplace-Stieltjes transform.
\begin{thm}\label{Ramamodify}
{\sl Let $X$ be a random variable with a symmetric distribution and finite moments $E|X|^n<\infty$ of all orders, $n\ge 0$. Then,
\begin{equation}\label{ramam}
  \int\limits_0^{\infty} t^{s-1}e^{-at}\phi_X(t)dt=\Gamma(s)E(a+iX)^{-s},\quad\hbox{for }a>0, \hbox{ and complex }s\hbox{ with }Re(s)>0,
\end{equation}
where $\phi_X(t)=E(e^{itX})$ is the characteristic function of $X$.
In particular, when $s=n>0$, we have for any $a>0$,
\begin{equation}\label{ramam2}
  e^{-at}\phi_X(t)=\sum\limits_{n=0}^{\infty}\frac{(-t)^nE(a+iX)^n}{n!},\quad t\ge 0.
\end{equation}
 }
 \end{thm}
 \begin{rem}
 Note that when $t=0$, both sides of (\ref{ramam2}) are equal to one. Note also that the integral $\frac1{\Gamma(s)}\int\limits_0^{\infty} t^{s-1}e^{-at}\phi_X(t)dt=h_X(s,a)$ is the Hurwitz Zeta function induced by $X$.
 \end{rem}
 \begin{proof}
 This theorem can be proved by the following two approaches:
 \begin{description}
   \item[Method 1.] First, we note that
   \begin{eqnarray*}
     e^{-ax}\phi_X(x)&=& \int\limits_{\Re} e^{-(a-it)x}dF_X(x) \\
      &=& \int\limits_{\Re}\sum\limits_{n=0}^{\infty}\frac{(-x)^n}{n!}(a-it)^ndF_X(t)\\
      &=& \sum\limits_{n=0}^{\infty}\int\limits_{\Re}\frac{(-x)^n}{n!}(a-it)^ndF_X(t)\\
      &=& \sum\limits_{n=0}^{\infty}\frac{(-x)^n}{n!}E(a-iX)^n=\sum\limits_{n=0}^{\infty}\frac{(-x)^n}{n!}E(a+iX)^n.
   \end{eqnarray*}
   Applying Ramanujan's Master Theorem, we have
   \begin{eqnarray*}
     \int\limits_0^{\infty} x^{s-1}  e^{-ax}\phi_X(x)dx&=& \int\limits_0^{\infty} x^{s-1}\sum\limits_{n=0}^{\infty}\frac{(-x)^n}{n!}E(a+iX)^ndx \\
     &=&\Gamma(s)E(a+iX)^{-s},\quad a>0,\quad Re(s)>0.
   \end{eqnarray*}
   \item[Method 2.] Directly by the characteristic function of $X$, we have
   \begin{eqnarray*}
    \int\limits_0^{\infty} x^{s-1}  e^{-ax}\phi_X(x)dx&=& \int\limits_0^{\infty} x^{s-1}  e^{-ax}\int\limits_{\Re} e^{itx}dF_X(t)\,dx\\
    &=& \int\limits_{\Re}\int\limits_0^{\infty}x^{s-1}e^{-(a-it)x}dx\,dF_X(t)=\int\limits_{\Re}\frac{\Gamma(s)}{(a-it)^s}dF_X(t)\\
    &=& \Gamma(s)E(a-iX)^{-s}= \Gamma(s)E(a+iX)^{-s}.
   \end{eqnarray*}
 \end{description}
 Since $E|X|^n<\infty$ for all $n\ge 1$, by analytic continuation, this identity can be extended to the range as $a>0$ and $Re(s)>0$.
 \end{proof}
 We can also obtain the following result concerning the Laplace-Stieltjes transform of a distribution function $G$.
 \begin{cor}\label{laplace}
 {\sl Let $Y\ge 0$ have a distribution function $G$, finite moments $E(Y^n)$ of all orders, $n\ge 0$, and moment generating function
 $L(\lambda)=E(e^{-\lambda Y})$, for $\lambda\ge 0$. Then,
 \[\int\limits_0^{\infty}x^{s-1}e^{-ax}L(x)dx=\Gamma(s)E(a+Y)^{-s},\quad\hbox{for }a>0,Re(s)>0.\]
 In particular, when $s=n\ge 0$, we have for $a> 0$ and $x\ge 0$,
 \[e^{-ax}L(x)=\sum\limits_{n=0}^{\infty}\frac{(-x)^n}{n!}E(a+Y)^n.\]
 }
 \end{cor}
 \begin{proof}
 The corollary can be proved in two different ways:
   \begin{description}
     \item [Method 1.]
     \begin{eqnarray*}
       e^{-ax}L(x) &=& e^{-ax}\int\limits_0^{\infty}e^{-xy}dG(y)=\int\limits_0^{\infty}e^{-x(a+y)}dG(y)\\
        &=& \int\limits_0^{\infty} \sum\limits_{n=0}^{\infty} \frac{(-x)^n}{n!}(a+y)^ndG(y)\\
        &=& \sum\limits_{n=0}^{\infty} \frac{(-x)^n}{n!}\int\limits_0^{\infty}(a+y)^ndG(y).
     \end{eqnarray*}
     Applying the Master Theorem, we have
     \begin{equation}\label{m1}
        \int\limits_{0}^{\infty}x^{s-1}e^{-ax}L(x)dx=\Gamma(s)\int\limits_0^{\infty}(a+y)^{-s}dG(y).
     \end{equation}
     \item [Method 2.] Prove the corollary directly.
      \begin{eqnarray*}
        &&\int\limits_{0}^{\infty}x^{s-1}e^{-ax}L(x)dx=\int\limits_0^{\infty}x^{s-1}e^{-ax}\int\limits_0^{\infty}e^{-xy}dG(y)\,dx\\
        &=&\int\limits_{0}^{\infty}\int\limits_{0}^{\infty} x^{s-1}e^{-(a+y)x}dx\,dG(y)\\
        &=&\int\limits_{0}^{\infty} \frac{\Gamma(s)}{(a+y)^s}dG(s)=\Gamma(s)\int\limits_{0}^{\infty}(a+y)^{-s}dG(y),\quad Re(s)>0.
     \end{eqnarray*}
   \end{description}
 \end{proof}
 \begin{ex}
   For standard normal distribution function $\Phi(x)=\int\limits_{-\infty}^x \frac1{\sqrt{2\pi}}e^{-\frac{t^2}2}dt$, $x\in\Re$.
   Its corresponding characteristic function is $\phi(t)=e^{-t^2/2}$, $t\in \Re$. We have
   \begin{eqnarray}
     && e^{-at}e^{-\frac12 t^2} =\sum\limits_{n=0}^{\infty} \frac{(-t)^n}{n!}\int\limits_{\Re}(a+ix)^n d\Phi(x),\quad t>0,\hbox{ and } \\
      && \int\limits_{0}^{\infty} x^{s-1}e^{-ax}e^{-\frac12 x^2}dx=\Gamma(s)\int\limits_{\Re}(a+ix)^{-s} d\Phi(x),\quad Re(s)>0.
   \end{eqnarray}
 \end{ex}
 \begin{ex}
 {\sl Consider the p.d.f. of the Gamma distribution
 \begin{equation}\label{gamma}
   f_{\alpha}(x)=\frac1{\Gamma(\alpha)}x^{\alpha-1}e^{-x},\quad x>0,\alpha>0,
 \end{equation}
 and its probability Laplace-Stieltjes transform
 \begin{equation}\label{LSt}
   L_{\alpha}(u)=\Big(\frac1{1+u}\Big)^{\alpha},\quad u>0.
 \end{equation}
 We have, for $a>0$,
 \begin{eqnarray}
   \frac{e^{-ax}}{(1+x)^{\alpha}} &=& \sum\limits_{n=0}^{\infty}\frac{(-x)^n}{n!}\int\limits_0^{\infty}(a+x)^nf_{\alpha}(x)dx,\quad x>0\quad{and} \\
   \int\limits_0^{\infty}x^{s-1} \frac{e^{-ax}}{(1+x)^{\alpha}}dx &=& \Gamma(s)\int\limits_0^{\infty}(a+x)^{-s}f_{\alpha}(x)dx,\quad Re(s)>0.
 \end{eqnarray}
 Furthermore, we have the interesting interchangeable relation
 \begin{equation}\label{exchange}
   \int\limits_0^{\infty} L_{\alpha}(x)f_{\beta}(x)dx=\int\limits_0^{\infty}L_{\beta}(x)f_{\alpha}(x)dx,\quad \alpha,\beta>0.
 \end{equation}
 }
 \end{ex}
The equations (3.10) and (3.11) can be derived by Corollary \ref{laplace}. Now we prove the interchangeable relation (\ref{exchange}).
\begin{eqnarray*}
&& \int\limits_0^{\infty} L_{\alpha}(x)f_{\beta}(x)dx\\
&=&  \int\limits_0^{\infty}\Big(\frac1{1+x}\Big)^{\alpha}\frac1{\Gamma(\beta)}x^{\beta-1}e^{-x}dx\\
&=&\int\limits_0^{\infty}\int\limits_0^{\infty}e^{-xy}\frac1{\Gamma(\alpha)}y^{\alpha-1}e^{-y}dy\frac1{\Gamma(\beta)}x^{\beta-1}e^{-x}dx
\end{eqnarray*}
\begin{eqnarray*}
   &=& \frac1{\Gamma(\alpha)\Gamma(\beta)} \int\limits_0^{\infty}\int\limits_0^{\infty}x^{\beta-1}e^{-x}e^{-xy}dx\,y^{\alpha-1}e^{-y}dy\\
   &=&  \frac1{\Gamma(\alpha)\Gamma(\beta)} \int\limits_0^{\infty}\int\limits_0^{\infty}x^{\beta-1}e^{-(1+y)x}dx\,y^{\alpha-1}e^{-y}dy\\
   &=&  \frac1{\Gamma(\alpha)\Gamma(\beta)} \int\limits_0^{\infty} \frac{\Gamma(\beta)}{(1+y)^{\beta}}y^{\alpha-1}e^{-y}dy \\
   &=& \int\limits_0^{\infty} \Big(\frac1{1+y}\Big)^{\beta}\frac1{\Gamma(\alpha)}y^{\alpha-1}e^{-y}dy\\
   &=& \int\limits_0^{\infty}L_{\beta}(y)f_{\alpha}(y)dy.
\end{eqnarray*}
The classical integral representation of the Beta function is
\begin{eqnarray*}
  B(p,q) &=& \frac{\Gamma(p)\Gamma(q)}{\Gamma(p+q)}=\int\limits_0^1 u^{p-1}(1-u)^{q-1}du \\
   &=& \int\limits_{\Re} \frac{ce^{-qcy}}{(1+e^{-cy})^{p+q}}dy,
\end{eqnarray*}
where $c>0$, $Re(q)>0$, $Re(p)>0$. Pitman and Yor (2003, p.300) yielded the p.d.f. of the $\cosh$ $\hat{C}_{\alpha}$ distribution
\begin{equation}\label{coshpdf}
  \psi_{\alpha}(x)=\frac{2^{\alpha-2}}{\pi} B(\frac{\alpha+ix}{2},\frac{\alpha-ix}{2})=\frac{2^{\alpha-2}}{\Gamma(\alpha)\pi}\Big|\Gamma(\frac{\alpha+ix}2)\Big|^2.
\end{equation}
As mentioned in Pitman and Yor (2003), any explicit formula of the p.d.f. of the sinh $\hat{S}_{\alpha}$ distribution or the tanh $\hat{T}_{\alpha}$ distribution is still
not available. In the following example, we will illustrate how to derive the p.d.f. of the sinh $\hat{S}_{\alpha}$ distribution and the tanh $\hat{T}_{\alpha}$ distribution.
\begin{ex}
{\sl The p.d.f.'s of the hyperbolic cosh, sinh, and tanh distributions are given below:
\begin{description}
  \item[(1)] The cosh $\hat{C}_{\alpha}$ distribution:
  \begin{eqnarray}
    \int\limits_{\Re} e^{itx}\psi_{\alpha}(x)dx &=& \Big(\frac1{\cosh{t}}\Big)^{\alpha},\quad t\in\Re,\quad \alpha>0, \\
    \hbox{ and the p.d.f. is } \psi_{\alpha}(x)&=& \frac{2^{\alpha-2}}{\pi} B(\frac{\alpha+ix}{2},\frac{\alpha-ix}{2})\nonumber\\
    &=&\frac{2^{\alpha-2}}{\Gamma(\alpha)\pi}\Big|\Gamma(\frac{\alpha+ix}2)\Big|^2.
  \end{eqnarray}
  \item[(2)] The sinh $\hat{S}_{\alpha}$ distribution:
  \begin{eqnarray}
   && \int\limits_{\Re} e^{itx}\phi_{\alpha}(x)dx = \Big(\frac{t}{\sinh{t}}\Big)^{\alpha},\quad t\in\Re,\quad \alpha>0, \\
     &&\hbox{ and the p.d.f. is }\nonumber\\
      &&\phi_{\alpha}(x)= \frac{2^{\alpha}}{\pi}\int\limits_0^1 \frac{\cos(x\ln{y})y^{\alpha-1}(\ln{y})^{\alpha}}{(y^2-1)^{\alpha}}dy,\nonumber\\
     &=& \frac{2^{\alpha}}{\pi}\int\limits_0^{\pi/2} \cos(x\ln(\sec{\theta}))\cot^{\alpha-1}(\theta)[\csc{\theta}\ln(\sec\theta)]^{\alpha}d\theta,\quad x\in\Re,\alpha>0.
  \end{eqnarray}
  \item[(3)] The tanh $\hat{T}_{\alpha}$ distribution:
   \begin{eqnarray}
   && \int\limits_{\Re} e^{itx}\eta_{\alpha}(x)dx = \Big(\frac{\tanh{t}}{t}\Big)^{\alpha},\quad t\in\Re,\quad \alpha>0, \\
     &&\hbox{ and the p.d.f. is }\nonumber\\
      &&\eta_{\alpha}(x)= \frac1{\pi}\int\limits_0^1 \frac{\cos(x\ln{y})(y^2-1)^{\alpha}}{(\ln{y})^{\alpha}(y^2+1)^2}\frac1{y}dy,\quad x\in\Re,\alpha>0.
  \end{eqnarray}
\end{description}
}
For the hyperbolic sinh-distributed random variable $\hat{S}_{\alpha}$, observe that
\begin{eqnarray*}
   && \int\limits_{\Re} e^{-ity}\Big(\frac{y}{\sinh{y}}\Big)^{\alpha}dy\\
   &=&\int\limits_{\Re} e^{-ity}\frac{y^{\alpha}}{[(e^y-e^{-y})/2]^{\alpha}}dy,\quad\hbox{by setting }y=\ln{x}, \\
   &=& \int\limits_{0}^{\infty} x^{-it}\frac{(\ln{x})^{\alpha}}{(x-\frac1x)^{\alpha}/2^{\alpha}}\frac1x\,dx \\
   &=&2^{\alpha}\int\limits_{0}^{\infty}\frac{x^{-it-1+\alpha}(\ln{x})^{\alpha}}{(x^2-1)^{\alpha}}dx\\
   &=&2^{\alpha}\Big[\int\limits_{0}^1\frac{x^{-it-1+\alpha}(\ln{x})^{\alpha}}{(x^2-1)^{\alpha}} dx+\int\limits_{1}^{\infty}\frac{x^{-it-1+\alpha}(\ln{x})^{\alpha}}{(x^2-1)^{\alpha}}dx\Big]\\
   &=&2^{\alpha}\Big[\int\limits_{0}^1\frac{x^{-it-1+\alpha}(\ln{x})^{\alpha}}{(x^2-1)^{\alpha}} dx+\int\limits_{0}^{1}\frac{y^{it+1-\alpha}(-\ln{y})^{\alpha}}{(\frac1{y^2}-1)^{\alpha}y^2}dy\Big],\quad\hbox{(by setting }y=\frac1{x})
  \end{eqnarray*}
\begin{eqnarray*}
     &=&2^{\alpha}\Big[\int\limits_{0}^1\frac{x^{-it-1+\alpha}(\ln{x})^{\alpha}}{(x^2-1)^{\alpha}} dx+\int\limits_{0}^{1}\frac{y^{it}y^{\alpha-1}(\ln{y})^{\alpha}}{(y^2-1)^{\alpha}}dy\Big]\\
     &=&2^{\alpha} \int\limits_{0}^1 \frac{[x^{-it}+x^{it}]x^{\alpha-1}(\ln{x})^{\alpha}}{(x^2-1)^{\alpha}}dx  \\
     &=& 2^{\alpha} \int\limits_0^1\frac{2\cos(t\ln{x})x^{\alpha-1}(\ln{x})^{\alpha}}{(x^2-1)^{\alpha}}dx
     \end{eqnarray*}
The above calculations lead to
\begin{eqnarray*}
&&\frac1{2\pi}\int\limits_{\Re} e^{-ity}\Big(\frac{y}{\sinh{y}}\Big)^{\alpha}dy\\
&=& \frac{2^{\alpha}}{\pi} \int\limits_0^1\frac{\cos(t\ln{x})x^{\alpha-1}(\ln{x})^{\alpha}}{(x^2-1)^{\alpha}}dx\\
&\equiv & \phi_{\alpha}(t),\quad t\in\Re,\quad \alpha>0.
\end{eqnarray*}
By Fourier transform, we know that
\[\int\limits_{\Re}e^{itx}\phi_{\alpha}(x)dx=\Big(\frac{t}{\sinh{t}}\Big)^{\alpha},\quad t\in\Re,\quad\alpha>0.\]
\begin{eqnarray*}
  \phi_{\alpha}(x) &=& \frac{2^{\alpha}}{\pi} \int\limits_0^1\frac{\cos(x\ln{y})y^{\alpha-1}(\ln{y})^{\alpha}}{(y^2-1)^{\alpha}}dy,\quad\hbox{by setting }y=\cos{\theta} \\
   &=& \frac{2^{\alpha}}{\pi}\int\limits_0^{\pi/2}\frac{\cos(x\ln(\cos{\theta}))\cos^{\alpha-1}{\theta}(\ln{cos{\theta}})^{\alpha}}
   {(-1)^{\alpha}(1-\cos^2{\theta})^{\alpha}}\sin{\theta}d\theta\\
    &=&  \frac{2^{\alpha}}{\pi}\int\limits_0^{\pi/2} \cos(x\ln{\sec{\theta}})\cot^{\alpha-1}{\theta}[\csc{\theta}\ln{\sec{\theta}}]^{\alpha}d\theta.
\end{eqnarray*}
Apropos of the tanh-distributed  random variable $\hat{T}_{\alpha}$, note that
\begin{eqnarray*}
   && \int\limits_{\Re} e^{-ity}\Big(\frac{\tanh{y}}{y}\Big)^{\alpha}dy \\
   &=& \int\limits_{\Re}  e^{-ity} \frac{[e^y-e^{-y}]^{\alpha}}{y^{\alpha}[e^y+e^{-y}]^{\alpha}}dy,\quad (\hbox{setting }x=e^y) \\
   &=& \int\limits_0^{\infty} x^{-it} \frac{(x^2-1)^{\alpha}}{(x^2+1)^{\alpha}(\ln{x})^{\alpha}}\frac1{x}dx\\
   &=&  \int\limits_0^1 x^{-it} \frac{(x^2-1)^{\alpha}}{(x^2+1)^{\alpha}(\ln{x})^{\alpha}}\frac1{x}dx+ \int\limits_1^{\infty} x^{-it} \frac{(x^2-1)^{\alpha}}{(x^2+1)^{\alpha}(\ln{x})^{\alpha}}\frac1{x}dx.
\end{eqnarray*}
Setting $u=\frac1{x}$, it becomes
\begin{eqnarray*}
   && \int\limits_0^1 x^{-it} \frac{(x^2-1)^{\alpha}}{(x^2+1)^{\alpha}(\ln{x})^{\alpha}}\frac1{x}dx+ \int\limits_0^1
   \frac{u^{it}(\frac1{u^2}-1)^{\alpha}}{(\frac1{u^2}+1)^{\alpha}(-\ln{u})^{\alpha}u^{-1}}\frac1{u^2}du\\
   &=&\int\limits_0^1 x^{-it} \frac{(x^2-1)^{\alpha}}{(x^2+1)^{\alpha}(\ln{x})^{\alpha}}\frac1{x}dx+ \int\limits_0^1
   u^{it}\frac{(u^2-1)^{\alpha}}{(u^2+1)^{\alpha}(\ln{u})^{\alpha}u}du \\
   &=&\int\limits_0^1 [x^{-it}+x^{it}] \frac{(x^2-1)^{\alpha}}{(x^2+1)^{\alpha}(\ln{x})^{\alpha}}\frac1{x}dx\\
   &=& \int\limits_0^1 2\cos(t\ln{x})\frac{(x^2-1)^{\alpha}}{(x^2+1)^{\alpha}(\ln{x})^{\alpha}}\frac1{x}dx,
\end{eqnarray*}
which implies
\begin{eqnarray*}
   &&\frac1{2\pi}\int\limits_{\Re} e^{-ity}\Big(\frac{\tanh{y}}{y}\Big)^{\alpha}dy \\
   &=& \frac1{\pi} \int\limits_0^1 \cos(t\ln{x})\frac{(x^2-1)^{\alpha}}{(x^2+1)^{\alpha}(\ln{x})^{\alpha}}\frac1{x}dx\equiv \eta_{\alpha}(t),\quad t\in\Re,\alpha>0.
\end{eqnarray*}
\end{ex}
The following three theorems are about the integral representation of the moments of the hyperbolic sinh-,cosh-, and tanh- distributed random variables.
\begin{thm}
  {\sl The hyperbolic sinh-moment function $S^*_{\alpha}(s,a)$ is defined by
  \begin{equation}\label{sinhmofun}
    S^*_{\alpha}(s,a)=\sum\limits_{n=0}^{\infty}\frac{\binom{n+\alpha-1}{n}}{(n+a)^s},\quad s>[\alpha],a>0,
  \end{equation}
  where $[\alpha]$ stands for the greatest integer less than or equal to $\alpha$.
  Then, $S_{\alpha}^*(s,a)$ has the integral representation
  \begin{eqnarray}
    S_{\alpha}^*(s,a) &=& \frac{\Gamma(s-\alpha)}{\Gamma(s)}\int\limits_{\Re}\frac1{(a-\frac{\alpha}2+\frac12 iy)^{s-\alpha}}\phi_{\alpha}(y)dy \\
     &=& \frac{\Gamma(s-\alpha)}{\Gamma(s)} E(a-\frac{\alpha}2+\frac12 iX_{\alpha})^{\alpha-s},
  \end{eqnarray}
  where the domain of $s$ can be extended analytically to the complex plane with $Re(s)>[\alpha]$, $a>\frac{\alpha}2$, and $X_{\alpha}$ is a random variable with the p.d.f.
  \begin{equation}\label{sinhpdf}
    \phi_{\alpha}(x)=\frac{2^{\alpha}}{\pi}\int\limits_0^1 \frac{\cos(x\ln{y})y^{\alpha-1}(\ln{y})^{\alpha}}{(y^2-1)^{\alpha}}dy,\quad x\in\Re,\alpha>0.
  \end{equation}
  }
\end{thm}
\begin{proof}
  By the definition of $S^*_{\alpha}(s,a)$, we know that the series converges absolutely for $Re(s)>[\alpha]$. Interchanging the sum and the integral, we have
  \begin{eqnarray*}
    \Gamma(s)S^*_{\alpha}(s,a) &=& \sum\limits_{n=0}^{\infty}\binom{n+\alpha-1}{n} \frac{\Gamma(s)}{(n+a)^s}  \\
     &=& \sum\limits_{n=0}^{\infty}\binom{n+\alpha-1}{n}\int\limits_0^{\infty} x^{s-1}e^{-(n+a)x}dx \\
     &=& \int\limits_0^{\infty} x^{s-1}e^{-ax}\sum\limits_{n=0}^{\infty}\binom{n+\alpha-1}{n}(e^{-x})^n dx \\
     &=& \int\limits_0^{\infty} x^{s-1}e^{-ax} \Big(\frac1{1-e^{-x}}\Big)^{\alpha}dx.
  \end{eqnarray*}
 Applying the formula $\int\limits_{\Re} e^{itx}\phi_{\alpha}(x)dx=\Big(\frac{t}{\sinh{t}}\Big)^{\alpha}$, we have
 \begin{eqnarray*}
    && \Gamma(s)S_{\alpha}^*(s,a)=\int\limits_0^{\infty} x^{s-1}e^{-ax}(\frac1{1-e^{-x}})^{\alpha}dx \\
    &=&\int\limits_0^{\infty} x^{s-1}e^{-ax}(\frac{e^{x/2}}{e^{x/2}-e^{-x/2}})^{\alpha} dx\\
    &=& \int\limits_0^{\infty} x^{s-1}e^{-ax} e^{\frac{\alpha x}2}x^{-\alpha} (\frac{x/2}{\sinh{\frac{x}2}})^{\alpha} dx \\
    &=& \int\limits_0^{\infty} x^{s-\alpha-1}e^{-(a-\frac{\alpha}2)x}\int\limits_{\Re} e^{ixy/2}\phi_{\alpha}(y)dy\,dx
 \end{eqnarray*}
 \begin{eqnarray*}
    &=& \int\limits_{\Re}\int\limits_0^{\infty} x^{s-\alpha-1}e^{-(a-\frac{\alpha}2-\frac{iy}2)x}dx\,\phi_{\alpha}(y)dy \\
    &=& \int\limits_{\Re}\frac{\Gamma(s-\alpha)}{(a-\frac{\alpha}2-\frac12 iy)^{s-\alpha}}\phi_{\alpha}(y)dy \\
    &=& \Gamma(s-\alpha)\int\limits_{\Re}(a-\frac{\alpha}2-\frac12 iy)^{\alpha-s} \phi_{\alpha}(y)dy \\
    &=& \Gamma(s-\alpha)E(a-\frac{\alpha}2+\frac12 iX_{\alpha})^{\alpha-s}.
 \end{eqnarray*}
 Since for all $n\ge 0$, the moments $E|X_{\alpha}|^n<\infty$, this equation can be extended analytically to $\{(s,a)\in \mathbb{C}\times \Re:\quad Re(s)>[\alpha], a>\alpha/2\}.$
\end{proof}
The following corollary is a direct result of the above theorem.
\begin{cor}
  {\sl The hyperbolic sinh-moment function $S_m^*(s,a)$ is defined by
  \begin{equation*}
    S_m^*(s,a)=\sum\limits_{n=0}^{\infty}\frac{\binom{n+m-1}{n}}{(n+a)^s},\quad s>m,\,a>0\quad{for }\quad m=1,2,\cdots.
  \end{equation*}
  Then, the integral representation of $S_m^*(s,a)$
  \begin{equation*}
    S_m^*(s,a)=\frac1{(s-1)\cdots (s-m)}\int\limits_{\Re} \frac1{(a-\frac{m}2+\frac12 iy)^{s-m}}\phi_m(y)dy,
  \end{equation*}
  can be extended analytically to $s\in\mathbb{C}$ with $Re(s)\ne 1,2,\cdots,m$, $a>m/2$, where $\phi_m(x)$ is a p.d.f. on $\Re$,
  \begin{equation*}
    \phi_m(x)=\frac{2^m}{\pi}\int\limits_0^1 \frac{\cos(x\ln{y})y^{m-1}(\ln{y})^m}{(y^2-1)^m}dy,\quad x\in\Re.
  \end{equation*}
  }
\end{cor}
\begin{thm}
{\sl Suppose that the hyperbolic cosh-moment function $C_{\alpha}^*(s,a)$ is defined by
\begin{equation}\label{coshmo}
  C_{\alpha}^*(s,a)=\sum\limits_{n=0}^{\infty} \frac{\binom{-\alpha}{n}}{(n+a)^s},\quad \alpha>0,\quad s>[\alpha],\quad a>0.
\end{equation}
Then it has the following integral representation
\begin{equation}\label{coshint}
  C_{\alpha}^*(s,a)=\int\limits_{\Re} \frac1{(a-\frac{\alpha}2+i\frac{x}2)^s}\frac1{2^{\alpha}}\psi_{\alpha}(x)dx,\quad s\in\mathbb{C},\quad a>\alpha/2,
\end{equation}
where $\psi_{\alpha}(x)=\frac{2^{\alpha-2}}{\pi}B(\frac{\alpha+ix}2,\frac{\alpha-ix}2)$, $x\in\Re$.
}
\end{thm}
\begin{proof}
By the definition of $C_{\alpha}^*(s,a)$, we can see that the series is absolutely convergent for $Re(s)>[\alpha]$.
Interchanging the sum and the integration, we have
\begin{eqnarray*}
   && \Gamma(s)C_{\alpha}^*(s,a)=\sum\limits_{n=0}^{\infty} \binom{-\alpha}{n} \int\limits_0^{\infty} x^{s-1}e^{-(n+a)x}dx \\
   &=& \int\limits_0^{\infty} x^{s-1}e^{-ax}\sum\limits_{n=0}^{\infty} \binom{-\alpha}{n}(e^{-x})^ndx \\
   &=& \int\limits_0^{\infty} x^{s-1}e^{-ax}\Big(\frac1{1+e^{-x}}\Big)^{\alpha}dx.
  \end{eqnarray*}
Now, since for $\alpha>0$,
\begin{equation*}
  \int\limits_{\Re} e^{itx}\psi_{\alpha}(x)dx=\Big(\frac1{\cosh{t}}\Big)^{\alpha},
\end{equation*}
we have
\begin{eqnarray*}
   && \Gamma(s)C_{\alpha}^*(s,a)=\int\limits_0^{\infty} x^{s-1}e^{-ax}\Big(\frac{e^{x/2}}{e^{x/2}+e^{-x/2}}\Big)^{\alpha}dx \\
   &=& \int\limits_0^{\infty} x^{s-1}e^{-(a-\alpha/2)x}\frac1{2^{\alpha}}\Big(\frac1{\cosh{\frac{x}2}}\Big)^{\alpha}dx \\
   &=& \int\limits_0^{\infty} x^{s-1}e^{-(a-\alpha/2)x}\frac1{2^{\alpha}}\int\limits_{\Re} e^{ixy/2}\psi_{\alpha}(y)dy\,dx \\
   &=& \int\limits_{\Re} \int\limits_0^{\infty} x^{s-1}e^{-(a-\frac{\alpha}2-\frac12 iy)x}dx\,\frac1{2^{\alpha}}\psi_{\alpha}(y)dy\\
   &=& \int\limits_{\Re} \frac{\Gamma(s)}{(a-\frac{\alpha}2-\frac12 iy)^s}\frac1{2^{\alpha}}\psi_{\alpha}(y)dy.
\end{eqnarray*}
By the symmetry of the p.d.f. $\psi_{\alpha}(y)$ and $E|X_{\alpha}|^n<\infty$ for all $n\ge 0$, we obtain the integral representation  (\ref{coshint}) and are able to extend it to include $s\in \mathbb{C}$, $a>\frac{\alpha}2$, $\alpha>0$.
\end{proof}
\begin{thm}
  {\sl Suppose that the hyperbolic tanh-moment function $T_{\alpha}^*(s,a)$ is defined by
  \begin{equation}\label{tanhmo}
    T_{\alpha}^*(s,a)=2^{\alpha}\sum\limits_{n=0}^{\infty}(-1)^n \binom{n+\alpha-1}{n}E(a+n+V_{\alpha})^{-s},
  \end{equation}
  where $\alpha$ is a positive integer, $s>\alpha$ is a real number,  $a>0$, and $V_{\alpha}$ is a random variable with
  \begin{equation*}
    V_{\alpha}\buildrel d\over{=}U_1+U_2+\cdots+U_{\alpha},
  \end{equation*}
  for which $U_1,\cdots,U_{\alpha}\buildrel {i.i.d.}\over{\sim} U(0,1)$.

  Then $T_{\alpha}^*(s,a)$ has the extended integral representation
  \begin{equation}\label{tanhint}
    T_{\alpha}^*(s,a)=E(a+\frac12 i X_{\alpha})^{-s},\quad s\in\mathbb{C},\quad a>0,
  \end{equation}
  where $X_{\alpha}$ is a random variable with the p.d.f.
  \begin{equation}\label{tanhpdf}
    \eta_{\alpha}(x)=\frac1{\pi}\int\limits_0^1 \frac{\cos(x\ln{y})(y^2-1)^{\alpha}}{(\ln{y})^{\alpha}(y^2+1)^2}\frac1{y}dy,\quad x\in\Re.
  \end{equation}
  }
\end{thm}
\begin{proof}
From the definition of the tanh-moment function $T_{\alpha}^*(s,a)$, we can see that it converges absolutely for those complex numbers $s$ with $Re(s)>[\alpha]$.
By interchanging the sum and the integral, we have
\begin{eqnarray*}
   && \Gamma(s)T_{\alpha}^*(s,a)=\Gamma(s)2^{\alpha}\sum\limits_{n=0}^{\infty}(-1)^n \binom{n+\alpha-1}{n}E(a+n+V_{\alpha})^{-s} \\
   &=&  \sum\limits_{n=0}^{\infty}2^{\alpha}(-1)^n \binom{n+\alpha-1}{n}\int\limits_{\Re}\int\limits_0^{\infty} x^{s-1}e^{-(a+n+v)x}dx\,dF_{V_{\alpha}}(v)\\
   &=& \int\limits_0^{\infty} x^{s-1}e^{-ax}\sum\limits_{n=0}^{\infty} 2^{\alpha}(-1)^n \binom{n+\alpha-1}{n}e^{-nx}\int\limits_{\Re}e^{-vx}dF_{V_{\alpha}}(v)\,dx \\
   &=& \int\limits_0^{\infty} x^{s-1}e^{-ax}2^{\alpha}\sum\limits_{n=0}^{\infty}\binom{-\alpha}{n} (e^{-x})^n\Big(\frac{1-e^{-x}}{x}\Big)^{\alpha}dx\\
   &=& \int\limits_0^{\infty} x^{s-1}e^{-ax}2^{\alpha}(1+e^{-x})^{-\alpha} \Big(\frac{1-e^{-x}}{x}\Big)^{\alpha}dx \\
   &=& \int\limits_0^{\infty} x^{s-1}e^{-ax} (\frac2{x})^{\alpha} \Big(\frac{1-e^{-x}}{1+e^{-x}}\Big)^{\alpha}dx \\
   &=& \int\limits_0^{\infty} x^{s-1}e^{-ax} \Big(\frac{\tanh{\frac{x}2}}{x/2}\Big)^{\alpha}dx.
\end{eqnarray*}
By the characteristic function of $\eta_{\alpha}$, we have
\begin{equation*}
  \int\limits_{\Re} e^{itx}\eta_{\alpha}(x)dx=\Big(\frac{\tanh{t}}{t}\Big)^{\alpha},\quad\hbox{and }
\end{equation*}
\begin{eqnarray*}
   &&\Gamma(s)T_{\alpha}^*(s,a)=\int\limits_0^{\infty} x^{s-1}e^{-ax} \Big(\frac{\tanh{\frac{x}2}}{x/2}\Big)^{\alpha}dx \\
   &=& \int\limits_0^{\infty} x^{s-1}e^{-ax} \int\limits_{\Re} e^{ixy/2}\eta_{\alpha}(y)dy\,dx \\
    &=&  \int\limits_0^{\infty}\int\limits_{\Re}x^{s-1}e^{-ax}e^{ixy/2}\eta_{\alpha}(y)dy\,dx \\
   &=& \int\limits_{\Re}\int\limits_0^{\infty}x^{s-1}e^{-(a-\frac12 iy)x}dx\,\eta_{\alpha}(y)dy\\
   &=&\int\limits_{\Re}\frac{\Gamma(s)}{(a-\frac12 iy)^s}\eta_{\alpha}(y)dy\\
   &=&\Gamma(s)E(a+\frac12i X_{\alpha})^{-s}.
 \end{eqnarray*}
   The last equality is due to the symmetry of the distribution of $X_{\alpha}$.
Since all moments $E|X_{\alpha}|^n<\infty$, $n\ge 0$, the integral representation (\ref{tanhint}) can be analytically extended to include
$s\in \mathbb{C}$, $a>0$.
\end{proof}

\end{document}